\documentclass[a4paper,10pt]{article}

\usepackage[english]{babel} %% Русские автопереносы и стандарты
\usepackage{amssymb,amsmath} %% Подключение AmsLaTeX не помешает
 \usepackage{graphicx}

\def\qed{\hfill $\square$}

\newtheorem{theorem}{Theorem}
\newtheorem{lemma}{Lemma}

\newcommand{\eps}{\varepsilon}

\newcommand{\p}{\mathcal{P}}

\title{Asymptotic stability of forced oscillations\\ emanating
from a limit cycle}
\author{{Oleg Makarenkov$^1$\thanks{supported by
the Grant BF6M10 of Russian Federation Ministry of Education and
U.S. CRDF (BRHE), by RFBR Grant 09-01-00468, by the President of
Russian Federation Young Researcher grant MK-1620.2008.1 and by
Marie Curie grant PIIF-GA-2008-221331.} and Rafael
Ortega$^2$\thanks{supported by D.G.I. MTM2008-02502,
Ministerio de Educaci\'on y Cultura, Spain }}\\[6mm]
{\small $^1$Department of Mathematics, Imperial College}\\
{\small London SW7 2AZ, London, United Kingdom}
\\
{\small {\tt  o.makarenkov@imperial.ac.uk}} \\ {\small
$^2$Departamento de Matem\'atica Aplicada, Facultad de Ciencias}
\\
{\small Universidad de Granada, 18071-Granada, Spain.}
\\
{\small {\tt   rortega@ugr.es }}} 
\date{}

\begin{document}
\maketitle

\section{Introduction}
The study of forced oscillations emanating from a limit cycle is a
classical problem in the theory of bifurcation. Around 1950 the
basic method to deal with this problem was developed by Malkin in
\cite{malb} and this study was continued by Loud in \cite{loud}.
The state of the art before the contributions of Malkin and Loud
can be found in the book by Lefschetz \cite{Lef}. To describe the
general framework we start with an autonomous system
$$\dot{x}=f(x)$$ having a closed orbit $\Gamma$ associated to a
periodic solution $x_0 (t)$ with period $T>0$. Notice that $T$ is
not necessarily the minimal period. The perturbation considered is
$$\dot{x}=f(x)+\eps g(t,x;\eps )$$ where $g$ is periodic in $t$
and its period is precisely $T$. The beginning of Malkin's method
is the construction of a $T$-periodic function $M=M(\theta )$
depending upon $x_0 (t)$ and $g(\cdot ,\cdot ;0)$. The zeros of
$M$ are intimately linked to the possible bifurcations to
$T$-periodic solutions for $\eps
>0$. Assuming some non-degeneracy conditions on $x_0 (t)$ one can
prove that if $\theta_*$ is a non-degenerate zero of $M$
($M(\theta_* )=0$, $M'(\theta_* )\neq 0$) then the perturbed
system has a family of $T$-periodic solutions satisfying
$$x_{\eps} (t)=x_0 (t+\theta_* )+O(\eps ),\; \; {\rm as}\; \eps
\downarrow 0.$$ It is also possible to analyze the case of a zero
of higher multiplicity ($M(\theta_* )=0$, $M'(\theta_* )=0,\cdots
,M^{(k-1)} (\theta_* )=0$, $M^{(k)} (\theta_* )\neq 0$) but this
requires long computations, see e.g. \cite{loud} and \cite
{kopnin}. More recently a topological approach has been taken in
\cite{Oleg}. A bifurcation exists as soon as $\theta_*$ is a zero
where $M$ changes sign. The next step after the existence of
bifurcating branches is the study of the stability properties.
This was already considered in \cite{malb}, \cite{loud} and
\cite{kopnin}. Assuming that $\Gamma$ is an exponential attractor
it can be proved that the bifurcating periodic solution is
asymptotically stable when $M'(\theta_* )>0$ and unstable when
$M'(\theta_* )<0$. If $\theta_*$ is a zero of a higher
multiplicity, then the implicit function approach taken in
\cite{loud} and \cite{kopnin} does not allow to detect bifurcation
of stable periodic solution on the basis of the sign of
$M^{(k)}(\theta_*)$ and some further computations have to be done.
See in particular equations (3.5) in \cite{kopnin} and (4.23) in
\cite{loud}. The purpose of our paper is to obtain a topological
version of this result for increasing or decreasing zeros when the
derivative of $M$ at $\theta_*$ can vanish. In particular, we are
interested in an unified answer which does not depend on the
multiplicity of $\theta_*$. We will get a positive answer in the
case of analytic systems. For this class of systems we will use a
variant of Lyapunov-Schmidt reduction that will allow us to prove
that if $M$ is not identically zero then the number of
$T$-periodic solutions is finite. This is inspired by the results
of Nakajima and Seifert in \cite{ns} and R.A. Smith in \cite{ras}.
Once we know that $T$-periodic solutions are isolated we can talk
about their topological index. This is just a localized version of
the topological degree and the connections of this index with the
stability properties of the corresponding solutions have been
discussed in \cite{Trans, krazab, kolesov, Tampa}. The computation
of the index is then obtained via a result in the line of those in
\cite{Oleg}.\par The rest of the paper is organized in three
sections. In Section \ref{dos} we present some preliminary results
on the autonomous system. The main Theorem as well as an example
illustrating its applicability can be found in Section~\ref{tres}.
This section also shows how to prove the main result via
topological degree. Finally Section~\ref{cuatro} is devoted to the
proofs of three Lemmas previously employed.

\section{The autonomous system}\label{dos}

In this section we present some elementary facts about the
non-perturbed system. They will be needed later in order to state
our main Theorem. Let us start with the autonomous system
\begin{equation}\label{aut}
  \dot x=f(x)
\end{equation}
defined on an open subset $\Omega$ of $\mathbb{R}^n.$ The vector
field $f:\Omega\to\mathbb{R}^n$ is real analytic.
\par
Assume that $x_0(t)$ is a non-constant periodic solution of
(\ref{aut}) with period $T>0.$ The associated variational equation
is
\begin{equation}\label{var}
  \dot y=f'(x_0(t))y.
\end{equation}
This is a $T$-periodic equation having the solution $\dot x_0(t).$
The Floquet multipliers are labelled as $\mu_1,...,\mu_n$ and
counted according to their multiplicity. It will be assumed that
they satisfy
\begin{equation}\label{floquet}
\mu_1=1,\ |\mu_2|<1,\ ...,|\mu_n|<1.
\end{equation}
This condition implies that the closed orbit $
  \Gamma=\{x_0(t),\ t\in[0,T]\}
$ is an attractor (see \cite{CL}). The region of attraction is an
open neighborhood of $\Gamma$ which will be denoted by
$\mathcal{A}\subset\Omega.$
\par
In view of the condition on the Floquet multipliers we know that
the space of $T$-periodic solutions of (\ref{var}) has dimension
one. The same property must hold for the adjoint system
\begin{equation}\label{adj}
  \dot z =-f'(x_0(t))^*z.
\end{equation}
The next result will provide an orientation in the space of
$T$-periodic solutions of (\ref{adj}).
\begin{lemma}\label{signadj}
There exists an unique $T$-periodic solution
   $z_0(t)$ of (\ref{adj})
   satisfying
$$
  \left<\dot x_0(t),z_0(t)\right>=1,\quad{\rm for \ any\
  }t\in\mathbb{R}.
$$
\end{lemma}

\noindent {\bf Proof.}
 It is based on Perron's lemma \cite{perron} (see also \cite{dem}, Sec. III, \S 12).
 This result says that if $y(t)$ and
$z(t)$ are arbitrary solutions of (\ref{var}) and (\ref{adj}) then
\begin{equation}\nonumber
  \left<y(t),z(t)\right> \equiv \mbox{constant}.
\end{equation}
We will prove that if $z_1(t)$ is a non-trivial $T$-periodic
solution of (\ref{adj}) then
\begin{equation}\label{discussion}
  \left<\dot{x}(0),z_1(0)\right> \neq 0.
\end{equation}
Since the space of $T$-periodic solution has dimension one this
will complete the proof.

 To prove (\ref{discussion}) we find a $n
\times n$ matrix $S$ such that
\begin{equation}\nonumber
  S^{-1}Y(T)S=\left(\begin{array}{cccc}
    1 & 0 & \hdots & 0 \\
    \cline{2-4}
    0 & \multicolumn{1}{|c}{} & &  \\
    \vdots & \multicolumn{1}{|c}{} & A &  \\
    0 & \multicolumn{1}{|c}{} & &  \end{array}\right),
\end{equation}
where $Y(t)$ is the matrix solution of (\ref{var}) with $Y(0)=I_N$
and $\det (A-I) \neq 0$. From the definition of $S$ we have that
its first column $S_1$ is an eigenvector of $Y(T)$ corresponding
to the eigenvalue $\mu_1=1$. In particular $S_1$ is parallel to
$\dot{x}(0)$. Consider the matrix $\Sigma = (S_2|..|S_n)$ composed
by the remaining columns of $S$. From the definition of $S$ and
$A$,
\begin{equation}\nonumber
  Y(T)\Sigma=\Sigma A.
\end{equation}
Next we apply Perron's Lemma to the solutions $Y(t)S_i$ and
$z_1(t)$,
\begin{equation}\nonumber
  \left<z_1(0),S_i\right>=\left<z_1(0),Y(T)S_i\right>, \quad i=2,..,n.
\end{equation}
This implies
\begin{equation}\nonumber
  z_1(0)^* \Sigma = z_1(0)^* Y(T)\Sigma = z_1(0)^* \Sigma A.
\end{equation}
Hence $z_1(0)^*\Sigma(I-A)=0$ and so $z_1^*(0)\Sigma=0$. Now we
can conclude that (\ref{discussion}) holds, for otherwise we
should have $z_1^*(0)S=0$ which is impossible if $z_1(t)$ is
non-trivial.

\qed

As a simple example we consider the planar system
$$
 \dot x=\left(1-|x|^2\right)x+i|x|^2x,\quad
 x=x_1+ix_2\in\mathbb{C}.
$$
It has the periodic solution $x_0(t)={\rm e}^{it}$ whose orbit
$\Gamma=\mathbb{S}^1$ attracts $\mathcal{A}=\mathbb{C}-\{0\}.$ The
period is $T=2N\pi,$ where $N\ge 1$ is an integer arbitrarily
chosen. The variational equation along $x_0(t)$ is
$$
  \dot y=(-1+2i)y+(-1+i){\rm e}^{2it}\overline{y}
$$
and has the Floquet solutions
$$
  y_1(t)=\dot x_0(t)=i{\rm e}^{it},\ \ y_2(t)={\rm
  e}^{(-2+i)t}(-1+i).
$$
In consequence $\mu_1=1$ and $\mu_2={\rm e}^{-2T}.$ The
computation of $z_0(t)$ follows from the proof of
Lemma~\ref{signadj}. We know that
$$
  \left<y_1(t),z_0(t)\right>=1,\ \ \left<y_2(t),z_0(t)\right>={\rm
  constant}=0.
$$
The periodicity of ${\rm e}^{2t}y_2(t)$ and $z_0(t)$ implies that
this last constant must vanish. From these equations one obtains
that
$$
  z_0(t)=(1+i){\rm e}^{it}.
$$

\section{Main result and an example}\label{tres}

Let us consider the perturbed system
\begin{equation}\label{pert}
  \dot x=f(x)+\eps g(t,x,\eps),
\end{equation}
where
$g:\mathbb{R}\times\Omega\times[0,\eps_*]\mapsto\mathbb{R}^n$ is
continuous and $T$-periodic in $t.$ We also assume that for each
$t\in \mathbb{R}$ the function $g(t,\cdot ,\cdot )$ has partial
derivatives up to the second order with respect to $(x,\eps )$ and
these derivatives are continuous as functions of the three
variables $ (t,x,\eps)$. The most important assumption on the
regularity of $g$ will be the analyticity with respect to $x.$
This means that for each $x_*\in\Omega$ there exists $r>0$ such
that if $\|x-x_*\|<r$ then for $j=1,...,n$
$$
  g_j(t,x,\eps)=\sum_{\alpha\in\mathbb{N}^n}g_{\alpha,j}(t,\eps)(x-x_*)^\alpha,\quad
  t\in\mathbb{R},\ \eps\in[0,\eps_*].
$$
Here $\alpha=(\alpha_1,...,\alpha_n)$ is a multi-index and we
employ the notation for powers
$x^\alpha=x_1^{\alpha_1}\cdot...\cdot x_n^{\alpha_n}.$ The
coefficients $g_{\alpha,j}$ are continuous and $T$-periodic in $t$
and the convergence in the above series is uniform in $t$ and
$\eps.$ As in the previous Section the vector field $f$ is real
analytic on $\Omega$ and this is enough to guarantee that the
solutions of (\ref{pert}) depend analytically upon initial
conditions once $\eps$ and $t$ have been fixed (see \cite{Lef}).

Again $x_0(t)$ is a non-constant $T$-periodic solution of
(\ref{aut}) satisfying (\ref{floquet}). We consider the function
$$
  M(\theta)=\int_0^T
  \left<g(t,x_0(t+\theta),0),z_0(t+\theta)\right>dt,
$$
where $z_0$ is given by Lemma~\ref{signadj}.
 This function is $T$-periodic and
real analytic and so it will have a finite number of zeros in
$[0,T[$ unless it is identically zero.

Given $\theta_*\in[0,T[$ a zero of $M,$ $M(\theta_*)=0,$ we say
that ${\rm index}(M,\theta_*)=1$ if $M(\theta)\cdot
(\theta-\theta_*)>0$ when $\theta\not=\theta_*$ is close to
$\theta_*.$ When the inequality is reversed we say that ${\rm
index}(M,\theta_*)=-1.$ In any other case we say that ${\rm
index}(M,\theta_*)=0.$

\begin{theorem} \label{main} In the previous setting assume that
$M$ is not identically zero and let $U$ be a bounded and open set
satisfying
$$\Gamma\subset U\subset\overline{U}\subset\mathcal{A}$$
(Recall that $\Gamma$ is the closed orbit associated to $x_0(t)$
and $\mathcal{A}$ is its region of attraction). Then there exists
$\eps_0>0$ such that if $0<\eps\le \eps_0$ the system (\ref{pert})
has a finite number of $T$-periodic solutions passing through
$\overline{U}.$ Moreover, if $\theta_*$ is a zero of $M$ with
${\rm index}(M,\theta_*)\not=0$ then there exists a $T$-periodic
solution $x_\eps(t)$ of (\ref{pert}) with
$$
  x_\eps(t)-x_0(t+\theta_*)\to 0\quad{\rm as}\quad \eps\downarrow
  0,
$$
uniformly in $t\in\mathbb{R}.$ This solution is asymptotically
stable if ${\rm index}(M,\theta_*)=1$ and unstable if ${\rm
index}(M,\theta_*)=-1.$
\end{theorem}

To illustrate the result we consider the planar system
\begin{equation}\label{ex}
  \dot x=(1-|x|^2)x+i|x|^2 x+\eps(a(t)+b(t)x+c(t)\overline{x}),
\end{equation}
where $x\in\mathbb{C}$ and $a,b,c:\mathbb{R}\to\mathbb{C}$ are
continuous and $2\pi$-periodic. The autonomous system  ($\eps=0$)
was already analyzed in the previous section and we can now
construct the function $M$ for $x_0(t)={\rm e}^{it}$, $z_0
(t)=(1+i)e^{it}$ and $T=2\pi .$ A direct computation leads to the
formula
$$
  M(\theta )={\rm Re} [\int_0^{2\pi} (a(t)+b(t)e^{i(t+\theta )}+c(t)e^{-i(t+\theta
  )})
  (1-i)e^{-i(t+\theta )} dt]=
  2\pi {\rm Re}\left[\left(\widehat{a}_1{\rm
  e}^{-i\theta}+\widehat{b}_0+\widehat{c}_2{\rm
  e}^{-2i\theta}\right)(1-i)\right],
$$
where $\widehat{a}_{m},$ $\widehat{b}_{m}$ and $\widehat{c}_{m}$
refer to the Fourier coefficients  of $a,$ $b$ and $c$,  namely
$$
   \widehat{f}_{m}=\frac{1}{2\pi}\int_0^{2\pi}f(t){\rm
   e}^{-imt}dt.
$$
In principle Theorem~\ref{main} would provide information on a
bounded region $U$ whose closure is contained in
$\mathbb{C}-\{0\}.$ However the specific properties of (\ref{ex})
will allow us to deduce global results. To illustrate this we
first claim that for $0\leq \eps <1$ any $2\pi$-periodic solution
$x(t)$ will satisfy $$\max_{t\in \mathbb{R}} ||x(t)|| \leq \rho_+
:=[1+||a||_{\infty} +||b||_{\infty} +||c||_{\infty} ]^{1/2} .$$
Indeed if $t_*$ is an instant when $m:=\max ||x(t)||=||x(t_* )||$
then the derivative $\frac{d}{dt}||x(t)||^2 =2\left<
x(t),\dot{x}(t)\right>$ must vanish at $t_*$. From the equation
(\ref{ex}) we deduce that $$||x(t_* )||^4 =||x(t_* )||^2 +\eps
\left< a(t_* )+b(t_* )x(t_* )+c(t_* )\overline{x(t_*)},x(t_*
)\right> .$$ It is not restrictive to assume that $m>1$ and and by
dividing the latter equality by $m^2$ the claimed estimate
follows. Next we observe that $x\equiv 0$ is a $2\pi$-periodic
solution for $\eps =0$. The variational equation is $\dot{y}=y$
with Floquet multipliers $\mu_1 =\mu_2 =e^{2\pi}$. A standard
perturbation result guarantees the existence of some $\rho_- \in
(0,1)$ such that, for small $\eps$, there is a unique
$2\pi$-periodic solution $z_{\eps}(t)$ satisfying $\max ||z_{\eps}
(t)||\leq \rho_-$. Moreover this solution is unstable since all
Floquet multipliers are greater than one. Now we apply
Theorem~\ref{main} in the region $$U=\{ x\in \mathbb{C}:\; \rho_-
<||x||<\rho_+ \} .$$ The function $M$ can be expressed as a
trigonometric polynomial of the type
$$M(\theta )=\beta +\alpha \cos (\theta +\phi )+\gamma \cos
2(\theta +\varphi ),$$ with $\beta =2\pi {\rm Re}[\widehat{b}_0
(1-i)]$, $2\pi \widehat{a}_1 (1-i)=\alpha e^{-i\phi}$, $2\pi
\widehat{c}_2 (1-i)=\gamma e^{-2i\varphi}$. Now it is clear that
$M$ is not identically zero if and only if
$$
  |\widehat{a}_1|+|{\rm Re}[\widehat{b}_0 (1-i)]|+|\widehat{c}_2|>0.
$$
In such a case (\ref{ex}) has a finite number of $2\pi$-periodic
solutions passing through $\overline{U}$, say $N$. From the above
discussions we conclude that also the number of $2\pi$-periodic
solutions on the whole plane is finite, namely $N+1$. When the
function $M$ does not vanish we obtain a uniqueness result:
$z_{\eps}$ is the unique $2\pi$-periodic solution. When $M$
changes sign we obtain at least two additional $2\pi$-periodic
solutions, one asymptotically stable and one unstable. Summing up,
we observe that in this example the function $M$ gives conditions
for the existence and stability that are rather sharp. Notice also
that the function $M$ can have zeros of the type $M(\theta_0
)=M'(\theta_0 )=M''(\theta_0 )=0$, $M'''(\theta_0 )\neq 0$ and
they lead to an asymptotically stable solution.
\par Before the proof of the Theorem we will state three lemmas
that will be proved in the next section. Our first preliminary
result goes back to \cite[page 387]{malb} and \cite{loud}. It
shows that the zeros of the function $M$ are relevant for the
location of $T$-periodic solutions. We shall say that a solution
$x(t)$ passes through a set $S\subset \mathbb{R}^2$ if $x(t)\in S$
for some real $t$.

\begin{lemma}\label{prel} Assume that $\eps_k\downarrow 0$ is a given
sequence and let $x_k(t)$ be a $T$-periodic solution of
(\ref{pert}) with $\eps=\eps_k$ and passing through
$\overline{U}.$ Then it is possible to extract a subsequence
$\{x_k(t)\}$ and a number $\theta_*\in[0,T[$ such that
$M(\theta_*)=0$ and
$$
  x_k(t)-x_0(t+\theta_*)\to 0\quad{\rm as}\quad k\to\infty
$$
uniformly in $t\in\mathbb{R}.$
\end{lemma}
For the next statements it will be convenient to employ the
Poincar\'e map $\mathcal{P}_{\eps}$ associated to (\ref{pert}).
Denoting by $x(t;\zeta ,\eps )$ the solution of (\ref{pert})
satisfying $x(0)=\zeta$, we notice that for small $\eps$ and
$\zeta \in \overline{U}$ this solution is well defined in $[0,T]$.
This is a consequence of the theorem on continuous dependence
since $\overline{U}$ is compact and for $\eps =0$ the solutions
starting at $\overline{U} \subset \mathcal{A}$ are globally
defined in the future. This observation allows us to define
$$\mathcal{P}_{\eps} :\overline{U}\to \mathbb{R}^n ,\; \; \; \zeta
\mapsto x(T;\zeta ,\eps ).$$ This map is analytic and its fixed
points are in a one-to-one correspondence with the $T$-periodic
solutions starting at $\overline{U}$.
\begin{lemma}\label{lyap_shm} Assume that
$\theta_0 \in \mathbb{R}$ is an isolated zero of $M$, then there
exist $\eps_0>0$ and $R>0$ such that for any $\eps \in (0,\eps_0)$
the Poincar\'e map $\p_\eps$ of (\ref{pert}) has at most a finite
number of fixed points in $B_R(x_0(\theta_0)).$
\end{lemma}
The third preliminary result will establish a link between the
index of the zeros of $M$ and the fixed point index of the
Poincar\'e map. Results of this type were already obtained in
\cite{Oleg} but we will present later an independent proof. The
Brouwer degree of a map $f$ on a domain $\Omega$ will be denoted
by ${\rm deg} (f,\Omega )$. It is assumed that $\Omega$ is open
and bounded and $f$ does not vanish on its boundary.
\begin{lemma}\label{gradito} Assume that $\theta_0$ is an isolated zero of $M$ and $ \mathcal{V}$ is an open
neighborhood of $x_0 (\theta_0 )$. Then there exist a number
$\eps_{\star} >0$ and a family of open sets $V_\eps \subset
\mathbb{R}^n$, $\eps \in (0,\eps_{\star} )$, satisfying $$ x_0
(\theta_0 )\in V_{\eps},\; \; \; V_{\eps} \subset \mathcal{V}$$
and such that
$$
  {\rm deg}(id-\mathcal{P}_\eps,V_\eps )=-{\rm index}(\theta_0,M),\quad {\rm whenever } \ \eps\in(0,\eps_\star).
$$
\end{lemma}We are now in the position to prove Theorem~\ref{main}.

\noindent {\bf Proof of Theorem~\ref{main}.} If the function $M$
does not vanish then (\ref{pert}) has no $T$-periodic solutions
passing through $\overline{U}$ when $\eps>0$ is small enough. This
is a consequence of Lemma~\ref{prel}. From now on we assume that
$M$ vanishes somewhere. Let $T^* >0$ be the minimal period of $x_0
(t)$, so that $T=kT^*$ for some $k=1,2,\dots$ The function $M$ has
period $T^*$ and the sequence of zeros of $M$ on $[0,T^*[$ is
denoted by
$$
  0\le\theta_1<\theta_2<...<\theta_m<T^*.
$$
Another consequence of Lemma~\ref{prel} is that for small $\eps$
any $T$-periodic solution of (\ref{pert}) passing through
$\overline{U}$ must remain close to the orbit $\Gamma$ for all
time. In particular we can assume that all $T$-periodic solutions
passing though $\overline{U}$ have an initial condition
corresponding to a fixed point of $\mathcal{P}_{\eps}$.
\par \noindent\underline{Step 1.} There exists $\eps_1 >0$ such that
if $\eps \in (0,\eps_1 )$ then $\mathcal{P}_{\eps}$ has a finite
number of fixed points.
\par \noindent
Once again we apply Lemma~\ref{prel} and restrict $\eps$ so that
all the fixed points are contained in some of the balls
$B_R(x_0(\theta_i)),$ $i=1,...,m,$ where $R$ is given by
Lemma~\ref{lyap_shm}. The union of these balls contains all the
fixed points of $\mathcal{P}_{\eps}$ and we know by
Lemma~\ref{lyap_shm} that they contain a finite number of fixed
points.\par We can also assume that $R$ has been chosen so that
these balls are pairwise disjoint. This will be employed later and
it is possible since $T^*$ is the minimal period and so the points
$x_0(\theta_i )$ and $x_0(\theta_j )$ are different whenever
$i\neq j$.

After this step we can define the index of a $T$-periodic solution
passing through $\overline{U}.$ Assume that $x(t)$ is such a
solution for some $\eps\in (0,\eps_1).$ We can find an open set
$\mathcal{W}\subset U$ such that $x(0)\in\mathcal{W}$ is the only
fixed point of $\p_\eps$ lying on $\overline{\mathcal{W}}.$ The
index of $x(t)$ is defined as
$$
  \gamma_T(x)={\rm deg}(id-\p_\eps,\mathcal{W}).
$$
In principle this index could take any integer value but the
condition (\ref{floquet}) implies that
\begin{equation}\label{k}
  \gamma_T(x)\in\{-1,0,1\}.
\end{equation}
This fact was already noticed by Krasnoselskii in \cite{Trans}. We
refer to \cite{Trans} or \cite{Tampa} for the proof.

\noindent\underline{Step 2.} If  $x(t)$ is a $T$-periodic solution
of (\ref{pert}) passing through $\overline{U},$ then $x(t)$ is
asymptotically stable if $\gamma_T(x)=1$ and unstable if
$\gamma_T(x)\not=1.$

The condition (\ref{floquet}) and the continuity of the Floquet
multipliers with respect to parameters imply the existence of a
positive number $\sigma >0$ such that if $B(t)$ is a $T$-periodic
and continuous matrix with $||B(t)||\leq \sigma$ for all $t$ then
the system
$$\dot{y}=(f'(x_0 (t))+B(t))y$$ has Floquet multipliers $\mu_1^*
,\cdots ,\mu_n^*$ with $\mu_1^*$ positive and dominant and
$|\mu_i^* |<1$ for $i=2,\cdots ,n$. After a time translation we
conclude that the same property holds for the more general class
of systems
$$\dot{y}=(f'(x_0 (t+\theta ))+B(t))y,\; \; \max ||B(t)||<\sigma ,\; \; B(t+T)=B(t).$$
For small $\eps$ any $T$-periodic solution passing through
$\overline{U}$ has a variational equation in this class and so the
Floquet multipliers have the structure described above. The
conclusion of Step~2 is a consequence of \cite{kolesov} and
\cite{Tampa}.

\noindent\underline{Step 3.} Assume that ${\rm
index}(M,\theta_i)\not=0.$ Then for any $\eps\in(0,\eps_2]$ the
equation (\ref{pert}) has a $T$-periodic solution $x$ with
$$
  x(0)\in B_R(x_0(\theta_i))\quad{\rm and}\quad \gamma_T(x)=-{\rm
  index}(M,\theta_i).
$$

 This is a consequence of Lemma \ref{gradito}. Indeed we can find an
 open set $V_{\eps}\subset B_R(x_0(\theta_i))$ with
$$
  {\rm deg}(id-\p_\eps,V_{\eps})=-{\rm index}(M,\theta_i)
$$
and the additivity of the degree implies that
$$
  {\rm deg}(id-\p_\eps,V_{\eps})=\sum\limits_{j=1}^m\gamma_T(x_j),
$$
where $x_1,...,x_m$ are the $T$-periodic solutions of (\ref{pert})
with $x_j(0)\in V_{\eps}.$ The conclusion follows from (\ref{k}).
Notice that the convergence of this periodic solution to $x_0
(t+\theta_i )$ as $\eps \to 0$ is a consequence of Lemma
\ref{prel} since the balls $B_R (x_0 (\theta_i ))$ are pairwise
disjoint. This completes the proof of the Theorem.

\section{Proofs of the Lemmas}\label{cuatro}

\noindent {\bf Proof of Lemma~\ref{prel}.} We present a proof for
completeness. Since $x_k$ passes through $\overline{U}$ one can
find $\tau_k\in[0,T]$ such that $x_k(\tau_k)\in\overline{U}.$
After extracting subsequences we can assume that
$$
   \tau_k\to\tau\quad{\rm and}\quad x_k(\tau_k)\to\zeta.
$$
Let $\widehat{x}(t)$ denotes the solution of (\ref{aut}) with
initial condition $\widehat{x}(\tau)=\zeta.$ Since $\zeta$ is a
point in the region of attraction $\mathcal{A}$ we know that
$\widehat{x}(t)$ is well defined in $[\tau,\infty[.$ By continuous
dependence we know that $x_k(t)$ converges to $\widehat{x}(t)$ and
the convergence is uniform on every compact interval where
$\widehat{x}(t)$ is well defined. In particular this applies to
$[\tau,\tau+T]$ and so $\widehat{x}(\tau)={\rm lim}x_k(\tau)={\rm
lim} x_k(\tau+T)=\widehat{x}(\tau+T).$ This implies that
$\widehat{x}(t)$ is a periodic solution of (\ref{aut}). Since
$\mathcal{A}$ is invariant for (\ref{aut}) and
$\widehat{x}(\tau)\in\mathcal{A}$ we deduce that the closed orbit
associated to $\widehat{x}$ must be contained in $\mathcal{A}.$
This implies that this orbit is precisely $\Gamma$ and so there
exists $\theta_*\in[0,T[$ such that
$\widehat{x}(t)=x_0(t+\theta_*).$ In particular $x_k(0)\to
x_0(\theta_*).$ It remains to prove that $M(\theta_*)=0.$ To this
end we consider the map
$$
   \Phi(\zeta,\eps)=\p_\eps(\zeta)-\zeta,\quad
   \zeta\in\overline{U},\ \eps\in[0,\eps_0].
$$
This is a $C^1$ map and the derivative $D\Phi(\zeta,\eps)$ is an
$n\times(n+1)$ matrix. We claim that the rank of
$D\Phi(x_0(\theta_*),0)$ is strictly less then $n.$ Otherwise the
equation $\Phi(\zeta,\eps)=0$ should describe a curve in a small
neighborhood of $(x_0(\theta_*),0).$ However the set $\Phi=0$
contains the curve $(x_0(\theta),0)$ and also the set of points
$(x_k(0),\eps_k)$ accumulating on $(x_0(\theta_*),0).$ Once we
know that ${\rm rank}D\Phi(x_0(\theta_*),0)<n,$ it remains to
prove that
$$
   {\rm rank}D\Phi(x_0(\theta),0)=n\quad{\rm if}\quad
   M(\theta)\not=0.
$$
The partial derivative with respect to $\xi$ is the $n\times n$
matrix
$$
  \partial_\zeta\Phi(x_0(\theta),0)=Y(T+\theta)Y(\theta)^{-1}-I_n,
$$
where $Y(t)$ is the matrix solution of (\ref{var}) with
$Y(0)=I_n.$ Again, the Fredholm alternative for linear
endomorphisms is applied to deduce that
$$
  {\rm Im}\partial_\zeta\Phi(x_0(\theta),0)=\left[{\rm
  Ker}\left(\left[Y(\theta)^*\right]^{-1}Y(\theta+T)^*-I_n\right)\right]^\bot.
$$
The kernel in the above formula corresponds to the initial
conditions at time $t=\theta$ of the $T$-periodic solutions of
(\ref{adj}). Hence it is spanned by $z_0(\theta)$ and so
$$
  {\rm
  Im}\partial_\zeta\Phi(x_0(\theta),0)=\left\{\eta\in\mathbb{R}^n:\eta\bot
  z_0(\theta)\right\}.
$$
By differentiability with respect to parameters, the function
$y(t)=\partial_\eps x(t,\zeta,\eps)$ with $\zeta=x_0(\theta),$
$\eps=0$ solves
$$
  \dot y=f'(x_0(t+\theta))y+g(t,x_0(t+\theta),0),\quad y(0)=0.
$$
A direct computation shows that
$$
  \frac{d}{dt}\left<y(t),z_0(t+\theta)\right>=\left<g(t,x_0(t+\theta),0),z_0(t+\theta)\right>
$$
and, integrating over the period,
$$
  \left<y(T),z_0(\theta)\right>=M(\theta).
$$
When $M(\theta)\not=0$ the vector $y(T)$ is not in the range of
$\partial_\zeta\Phi(x_0(\theta),0)$ and so
$$
  {\rm
  rank}\left(\partial_\zeta\Phi(x_0(\theta),0)\left|\partial_\eps\Phi(x_0(\theta),0)\right.\right)=(n-1)+1=n.
$$
\qed

\noindent {\bf Proof of Lemma~\ref{lyap_shm}.} It is based on a
variant of the Lyapunov-Schmidt reduction. We divide it in four
steps.

\noindent\underline{1. The change of variables.} The dominant
eigenvalue of $L=(\p_0)'(x_0(\theta_0))$ is $\mu_1=1$ with
eigenvector $\dot{x}_0(\theta_0)$. This eigenvalue is simple and
so we can find a linear projection $\pi$ in $\mathbb{R}^n$
satisfying
\begin{equation}\nonumber
 \pi^2=\pi, \quad \pi L=L\pi, \quad {\rm Ker}~\pi=\{\lambda \dot{x}(\theta_0) \mbox{;} \ \lambda \in
 \mathbb{R}\}.
\end{equation}
This is so-called spectral projection and the hyperplane $Y={\rm
Im} (id- L)$
is invariant under $L$.\\
Moreover,
\begin{equation}\label{eigen}
 \sigma(L_Y)=\{\mu_2,.., \mu_n\},
\end{equation}
where $L_Y:Y \to Y$ is the restriction of $L$ to $Y$. In the rest
of the proof $v$ denotes a generic vector lying in $Y$.

Consider the map
\begin{equation}\nonumber
 \Phi:(\theta,v) \in \mathbb{R} \times Y \mapsto x_0(\theta)+v \in
 \mathbb{R}^n.
\end{equation}
This is an analytic function with partial derivatives at
$(\theta_0,0)$,
\begin{equation}\nonumber
  \partial_\theta \Phi (\theta_0,0)=\dot{x}_0(\theta_0), \quad
  \partial_Y \Phi (\theta_0,v)=id_Y.
\end{equation}
The Inverse Function Theorem implies that $\Phi$ is a local
diffeomorphism mapping $(\theta_0,0)$ onto $x_0(\theta_0)$. In a
neighborhood of this point we reduce the search of fixed points of
$\p_\eps$ to the equation $\p_\eps \circ \Phi =\Phi$. More
precisely we consider the equation
\begin{equation}\label{nv}
 \p_\eps(x_0(\theta)+v)=x_0(\theta)+v, \quad
 |\theta-\theta_0|<\Delta, \quad \|v\|< \Delta,
\end{equation}
for some small $\Delta>0$. Notice that $\Phi$ is independent of
$\eps$ and so $\Delta$ is uniform in $\eps \geqslant 0$.\\

\noindent\underline{2. The auxiliary equation.} The equation
(\ref{nv}) can be interpreted as a system in the unknowns $\theta$
and $v$. As usual we apply $\pi$ and solve in $v$. This means that
we look at the implicit function problem
\begin{equation}\nonumber
  F(\theta,v \mbox{;}\eps):=\pi \p_\eps (x_0(\theta)+v)-\pi
  x_0(\theta)-v=0.
\end{equation}
This function maps $|\theta-\theta_0|<\Delta, \quad \|v\|<\Delta,
\quad \eps \in [0,\eps_*]$ into $Y$ and satisfies
\begin{equation}\nonumber
  \partial_v F(\theta_0,0\mbox{;}0)=L_Y-id_Y.
\end{equation}
From the condition (\ref{eigen}) we deduce that the implicit
Function Theorem is applicable and so we find $r>0$ and
$\alpha:[\theta_0-r,\theta_0+r]\times[0,r] \to Y$ such that
\begin{equation}\nonumber
  \pi \p_\eps(x_0(\theta)+\alpha (\theta,\eps))=\pi
  x_0(\theta)+\alpha (\theta,\eps).
\end{equation}
Moreover this is the only solution of $F(\theta,v;\eps)=0$ in some
ball $\|v\|<R$. The function $\alpha$ is of class $C^1$ and
analytic with respect to $\theta$. Due to the uniqueness of
$\alpha$ we have $\alpha(\theta,0)=0$ for any $\theta \in
[\theta_0-r,\theta_0+r]$, which can be combined with the
smoothness of $\alpha$ to find a number $\mu >0$ such that
\begin{equation}\label{dif}
  \|\alpha(\theta,\eps)\| \leqslant \eps \mu \quad \mbox{for any} \quad  \theta
\in [\theta_0-r,\theta_0+r], \quad \eps \in [0,r].
\end{equation}
In this process it can be necessary to reduce the size of $r$.

\noindent\underline{3. The bifurcation equation.} Assume that $x(t
\mbox{;}\Xi,\eps)$ is a $T$-periodic solution of (\ref{pert}) with
$\Xi$ close to $x_0(\theta_0)$ and $\eps$ small and positive. We
know from the previous steps that the initial condition can be
expressed as
\begin{equation}\nonumber
  \Xi =x_0(\Theta)+\alpha(\Theta,\eps)
\end{equation}
for some $\Theta \in [\theta_0-r,\theta_0+r]$. Our next task is to
show that $\Theta$ must be a zero of the function
\begin{equation}\nonumber
  M_\eps (\theta):= \int_0^T \left< b_\eps (t,\theta),z_0(t+\theta)\right>dt
\end{equation}
with
\begin{eqnarray}\nonumber
  b_\eps(t,\theta)&:=&g(t,x(t,\xi,\eps),\eps)-\frac{1}{\eps}[f(x(t,\xi,\eps))
  -f(x_0(t+\theta))-f'(x_0(t+\theta))\cdot
   (x(t,\xi,\eps)-x_0(t+\theta))] \nonumber
\end{eqnarray}
and
\begin{equation}\nonumber
  \xi = x_0(\theta)+\alpha(\theta,\eps).
\end{equation}
By construction $y(t)=x(t,\Xi,\eps)-x_0(t+\Theta)$ has to be a
$T$-periodic solution of the linear equation
\begin{equation}\nonumber
  \dot{y}=f'(x_0(t+\Theta))y+\eps b_\eps (t,\Theta).
\end{equation}
The Fredholm alternative implies that $\Theta$ is a zero of
$M_\eps$.

\noindent\underline{4. Conclusion: the role of analyticity.} In
view of the previous steps it is enough to show that the function
$M_\eps$
has a finite number of zeros in $[\theta_0-r,\theta_0+r]$ for small $\eps$.\\
Since $\alpha(\theta,0)=0$ we obtain by continuous dependence that
\begin{equation}\nonumber
  b_\eps (t,\theta) \to g(t,x_0(t+\theta),0) \quad \mbox{as} \quad \eps \to
  0
\end{equation}
uniformly in $t \in [0,T]$ and $\theta \in
[\theta_0-r,\theta_0+r]$. Indeed we also need to use that $f$ is
smooth and the estimate (\ref{dif}). This is required to prove
that the term related to $f$ goes to zero. Also the
differentiability with respect to initial conditions and
parameters plays a role here.

 The function $M_\eps$ converges to
$M$ as $\eps \to 0$ uniformly in $\theta \in
[\theta_0-r,\theta_0+r]$. We are assuming that $M$ is not
identically zero and so the same must happen to $M_\eps$ for small
$\eps$. Since $M_\eps$ is analytic we conclude that it has a
finite numbers of zeros in $[\theta_0-r,\theta_0+r]$. This is
valid for $\eps \in ]0,\eps_0[$ with $\eps_0>0$ sufficiently
small.

\qed

{\bf Remark} The standard Lyapunov-Schmidt reduction for the
equation $P_\eps(\xi)=\xi$ would start with the splitting
\begin{equation}\nonumber
  \xi=\eta \dot{x}_0(\theta_0)+v, \qquad \eta \in \mathbb{R}, \ v \in Y,
\end{equation}
and considering the system
\begin{equation}\nonumber
\left\{ \begin{array}{l} \pi \p_\eps (\eta
\dot{x}_0(\theta_0)+v)=v
\\ (id-\pi)\p_\eps (\eta
\dot{x}_0(\theta_0)+v)=\eta\dot{x}_0(\theta_0).  \end{array}
\right.
\end{equation}
Instead of this we are considering a sort of nonlinear splitting
induced by the change of variables of Step~1. The advantage is
that our bifurcation equation leads directly to $M(\theta)=0$ as
$\eps \downarrow 0$. The same approach is taken by Hale and Taboas
in \cite{Hale+Taboas}, but they prefer to work in an infinite
dimensional
framework.\\

\noindent{\bf Proof of Lemma~\ref{gradito}.} First we pick up any
$n-1$ linearly independent solutions $y_1,...,y_{n-1}$ of
(\ref{var}) whose initial conditions at $\theta_0$ satisfy $\left<
y_i(\theta_0),z_0(\theta_0)\right>=0$. Next we consider the
$n\times (n-1)$ matrix $Y_1(\theta)=(y_1(\theta)|\dots
|y_{n-1}(\theta))$ and notice that
\begin{equation}\label{ole}
  Y_1(\theta +T)=Y_1(\theta )A_{\theta}
  \end{equation}
  where $A_{\theta}$ is a $(n-1)\times (n-1)$ matrix with
  eigenvalues $\mu_2 ,\dots ,\mu_n$. To verify this it is enough
  to observe that the hyperplane $V_{\theta}$ spanned by $y_1
  (\theta ),\dots ,y_{n-1} (\theta )$ is invariant under the
  monodromy operator $M_{\theta}: \; y(\theta )\mapsto y(\theta
  +T)$. This is a consequence of Perron's Lemma. The eigenvector
  of $M_{\theta}$ associated to $\mu_1 =1$ is $\dot{x}_0 (\theta
  )$ and does not belong to $V_{\theta}$. In consequence the
  restriction of $M_{\theta}$ to $V_{\theta}$ has eigenvalues
  $\mu_2 ,\dots ,\mu_n$. The matrix $A_{\theta}$ is precisely the
  representation of this restriction with respect to the basis
  $y_1 (\theta ),\dots ,y_{n-1}(\theta )$. This property of
  the matrix $Y_1 (\theta )$ will be employed several times. First
  we will employ it to evaluate the topological degree of the auxiliary map
$$
  \Phi_\eps(\theta,\zeta)=-\eps M(\theta)\dot
  x_0(\theta)+(Y_1(\theta
  )-Y_1(\theta+T))\zeta
$$
with respect to $\Omega_\delta
:=(\theta_0-\delta,\theta_0+\delta)\times B_\delta(0),$ where
$\delta >0$ is a small number and $B_\delta(0)$ is an open ball in
$ \mathbb{R}^{n-1}.$ We will impose several restrictions on the
size of $\delta$, the first being that $M$ has no zeros on
$[\theta_0-\delta,\theta_0+\delta]$ other then $\theta_0$. Notice
that by linear independence the equation
$\Phi_\eps(\theta,\zeta)=0$ in $\overline{\Omega}_{\delta}$ splits
as $ M(\theta )=0$ and $(Y_1(\theta
  )-Y_1(\theta+T))\zeta =0$. Then $\theta =\theta_0$ and from the identity (\ref{ole}) we
  deduce that $\zeta =0$. Thus the degree we
  want to compute is well defined and does not change if we replace
  $\Omega_{\delta}$ by any sub-domain containing $(\theta_0 ,0)$. For the effective computation
  we diminish $\delta>0$ in such a way that $\Phi_\eps$ is linearly
homotopic to the vector field
$$
  \widehat{\Phi}(\theta,\zeta)=-M(\theta)\dot
  x_0(\theta_0)+(Y_1(\theta_0)-Y_1(\theta_0+T))\zeta
$$
for $\eps>0$ sufficiently small so that ${\rm
deg}(\Phi_\eps,\Omega_\delta)={\rm
deg}(\widehat{\Phi},\Omega_\delta).$ The matrix $S=(\dot{x}_0
(\theta_0 )| Y_1 (\theta_0 )-Y_1 (\theta_0 +T))$ is non-singular
and the map $\widehat{\Phi}$ can be expressed as $S\circ
[(-M)\times id]$. By the theorems on the evaluation of the
topological index of a composition of vector fields (see e.g.
\cite{krazab}, Theorem~7.1), of a product of vector fields (see
e.g. \cite{krazab}, Theorem~7.4) and of a linear vector field (see
\cite{krazab}, Theorem~6.1) we have that
\begin{eqnarray}
{\rm deg}(\widehat{\Phi},\Omega_\delta)= {\rm
index}((\theta_0,0),\widehat{\Phi})= {\rm index}(0,S) \cdot {\rm
index}((\theta_0,0),(-M)\times id)=  -{\rm sign}\det S \cdot{\rm
index}(\theta_0,M).\label{EQN}
\end{eqnarray}
Another restriction on $\delta$ that will be useful later is
related to the map $\psi (\theta ,\zeta )=x_0 (\theta )+Y_1
(\theta )\zeta$. This map must be a diffeomorphism from
$\Omega_{\delta}$ onto its image and $\psi
(\overline{\Omega}_{\delta} )\subset \mathcal{V}$. Notice that
this is possible since $\det \psi '(\theta_0 ,0)=\det (\dot{x}_0
(\theta_0 )|Y_1 (\theta_0 )\neq 0$.\par
 Our next step is to show that the vector fields
\begin{equation*}
  \mathcal{F}_{\eps} (\theta,\zeta)=(id-\mathcal{P}_\eps)(x_0(\theta)+Y_1(\theta)\zeta)
\end{equation*}
and $\Phi_\eps$ are homotopic on a sub-domain of $\Omega_\delta$
for $\eps>0$ sufficiently small. Let $x(t;\theta,\zeta,\eps)$ be
the solution of (\ref{pert}) satisfying
$x(0)=x_0(\theta)+Y_1(\theta)\zeta.$ The Taylor expansion leads to
\begin{equation*}
  x(t,\theta,\zeta,\eps)=x_0(t+\theta)+Y_1(t+\theta)\zeta+\eps\int_0^t
  Y(t+\theta)Y(s+\theta)^{-1}g(s,x_0(s+\theta),0)ds+O(\eps^2
  +\|\zeta\|^2),
\end{equation*}
where we recall that the matrix $Y(t)$ was defined in Section
\ref{dos}. This expansion is obtained by computing the derivatives
with respect to $\zeta$ and $\eps$ and applying the formula of
variation of constants. The matrix $Y^* (t+\theta )^{-1}Y^*
(\theta )$ is fundamental at $t=\theta$ for the adjoint system and
so
$$
  z_0(t+\theta)=Y^*(t+\theta)^{-1}Y^*(\theta)z_0(\theta).
$$
From the periodicity of $z_0$ we deduce that
$$
   Y^*(T+\theta)z_0(\theta)=Y^*(\theta)z_0(\theta).
$$
Thus,
\begin{eqnarray*}
   \left<\int_0^T
   Y(T+\theta)Y^{-1}(s+\theta)g(s,x_0(s+\theta),0)ds,z_0(\theta)\right>=\\
   =\int_0^T\left<g(s,x_0(s+\theta),0),z_0(s+\theta)\right>ds=M(\theta).
\end{eqnarray*}
In consequence,
$$
  \left< \mathcal{F}_{\eps} (\theta ,\zeta ), z_0(\theta)\right>=-\eps
  M(\theta)+O(\eps^2+\|\zeta\|^2),\; \; \;
 \mathcal{F}_{\eps} (\theta ,\zeta )=
   (Y_1(\theta)-Y_1(\theta+T))\zeta+\eps\gamma(\theta)+O(\eps^2
   +\|\zeta\|^2),
$$
where $\gamma$ is defined by an integral. Perhaps after a new
reduction of the size of $\delta$ we can find a positive constant
$\Lambda$ such that $$\max_{k=1,\dots ,n-1} |\left<
(Y_1(\theta)-Y_1(\theta+T))\zeta ,y_k (\theta )\right> |\geq
\Lambda ||\zeta ||, \; \; {\rm for\; every}\; \zeta \in
\mathbb{R}^{n-1}\; \; {\rm and}\; |\theta -\theta_0 |\leq \delta
.$$ To justify this assertion we notice that, by continuity, it is
enough to check it for $\theta =\theta_0$ and in this case it
follows from (\ref{ole}) since
$(Y_1(\theta_0)-Y_1(\theta_0+T))=Y_1 (\theta_0 )(I-A_{\theta_0})$
and $(I-A_{\theta_0})$ is non-singular. From now on the number
$\delta$ will be kept fixed. We are going to compute the degree of
$\mathcal{F}_{\eps}$ on the set $W_{\eps}=\{ (\theta ,\zeta ):\;
|\theta -\theta_0 |<\delta ,\; ||\zeta ||<\eps^{2/3}\}$. The
boundary of $W_{\eps}$ is composed by $\Delta_1: \theta =\theta_0
\pm \delta ,||\zeta ||\leq \eps^{2/3}$ and $\Delta_2: |\theta
-\theta_0 |\leq \delta ,||\zeta ||= \eps^{2/3}$. On $\Delta_1$ we
observe that for $\eps$ small enough
$${\rm sign} \left< \mathcal{F}_{\eps} (\theta ,\zeta ),z_0
(\theta )\right> =-{\rm sign} M(\theta ),\; \; {\rm with}\; \theta
=\theta_0 \pm \delta .$$ On $\Delta_2$ we claim that for some
$k=1,\dots ,n-1$ (depending on $\zeta$),
$${\rm sign} \left< \mathcal{F}_{\eps} (\theta ,\zeta ),y_k
(\theta )\right> =-{\rm sign} \left< (Y_1 (\theta )-Y_1 (\theta
+T)) \zeta ,y_k (\theta )\right> .$$ Indeed, from the expansion of
$\mathcal{F}_{\eps}$ we find that for each $k$
$$\left< \mathcal{F}_{\eps} (\theta ,\zeta ),y_k (\theta )\right>
=- \left< (Y_1 (\theta )-Y_1 (\theta +T)) \zeta ,y_k (\theta
)\right>+O(\eps ).$$ For some $k$, $ |\left< (Y_1 (\theta )-Y_1
(\theta +T)) \zeta ,y_k (\theta )\right> |\geq \Lambda \eps^{2/3}$
and this term is dominant, leading to the coincidence of the
signs. Summarizing, for $(\theta,\zeta)\in \partial W_{\eps}$ the
vectors $\Phi_\eps(\theta,\zeta)$ and
$\mathcal{F}_{\eps}(\theta,\zeta)$ do not point in opposite
directions and, therefore, the vector fields $\Phi_\eps$ and
$\mathcal{F}_{\eps}$ are linearly homotopic on $W_\eps$ (see
\cite[theorem~2.1]{krazab}). By excision,
\begin{equation}\label{FO}
  {\rm deg}(\mathcal{F}_{\eps},W_\eps)={\rm
  deg}(\widehat\Phi,\Omega_\delta)=-{\rm sign}\det S \cdot{\rm ind}(\theta_0,M).
\end{equation}
To finish the proof we define $V_{\eps}=\psi (W_{\eps})$ and
observe that $(id-\mathcal{P}_{\eps})\circ \psi
=\mathcal{F}_{\eps}$ on $W_{\eps}$. The theorem on the degree of
the composition implies that $${\rm deg} (id-\mathcal{P}_{\eps}
,V_{\eps})\cdot {\rm deg}(\psi -x_0 (\theta_0 ),W_{\eps})={\rm
deg}(\mathcal{F}_{\eps} ,W_{\eps}).$$ For instance, Theorem~7.2,
Formula~7.6 in \cite{krazab} is applicable since $\partial
V_{\eps} =\psi (\partial W_{\eps})$, $V_{\eps}$ is connected and
$x_0 (\theta_0 )\in V_{\eps}$. By the linearization theorem for
topological degree (see e.g. \cite{krazab}, Theorem~6.3) we have
that
\begin{equation}\label{H2}
  {\rm deg}(\psi -x_0 (\theta_0 ),W_\eps )={\rm sign}\det\psi'(\theta_0,0)={\rm
  sign}\det(\dot x_0(\theta_0)|Y_1(\theta_0)).
\end{equation}
The conclusion of the Lemma follows from these last identities and
(\ref{FO}) because
\begin{equation}\label{idsign} {\rm sign} \det (\dot{x}_0
(\theta_0 )|Y_1 (\theta_0 ))={\rm sign}\det S.
\end{equation}
To prove this claim we consider the family of matrices
$$Y_1 (\theta_0 )-\lambda Y_1 (\theta_0 +T)=Y_1 (\theta_0
)(I-\lambda A_{\theta_0} ),\; \; \; \lambda \in [0,1],$$ where
once again we have used (\ref{ole}). For $\lambda =0$ and $\lambda
=1$ we obtain the second blocks of the matrices appearing in the
identity (\ref{idsign}). The eigenvalues of $A_{\theta_0}$ are
$\mu_2 ,\dots ,\mu_n$, all of them with modulus less than one.
Hence $$\det (\dot{x}_0 (\theta_0 )|Y_1 (\theta_0 )-\lambda Y_1
(\theta_0 +T))\neq 0$$ for all $\lambda \in [0,1]$ and so the sign
of this determinant is independent of $\lambda$. The identity
(\ref{idsign}) expresses this fact for the extreme values of
$\lambda$. \qed

\end{document}